\documentclass[iicol,sn-mathphys]{sn-jnl}
\jyear{2023}
\usepackage{graphicx}
\graphicspath{{./Figs/}}
\usepackage{graphicx}%
\usepackage{multirow}%
\usepackage{amsmath,amssymb,amsfonts}%
\usepackage{amsthm}%
\usepackage{mathrsfs}%
\usepackage[figuresright]{rotating}%
\usepackage[title]{appendix}%
\usepackage{xcolor}%
\usepackage{textcomp}%
\usepackage{manyfoot}%
\usepackage{booktabs}%
\usepackage{dsfont}
\usepackage{amsmath}
\usepackage{enumitem}

\usepackage{tikz}
\usepackage{pgfplots}
\pgfplotsset{compat=newest}
\pgfplotsset{plot coordinates/math parser=true}
\pgfplotsset{grid style={dotted,gray}}
\usetikzlibrary{plotmarks,pgfplots.fillbetween}

\newcommand{\lc}{\left\{}
\newcommand{\rc}{\right\}}
\newcommand{\lp}{\left(}
\newcommand{\rp}{\right)}

\newcommand{\dd}{\,\mathrm{d}}
 
\theoremstyle{thmstyleone}%
\theoremstyle{thmstyletwo}%

\theoremstyle{thmstylethree}%

\begin{document}

\title[Article Title]{Topology Optimization of Broadband Acoustic Transition Section: A Comparison between Deterministic and Stochastic Approaches}

\author*[1]{\fnm{Abbas} \sur{Mousavi}}\email{sam@cs.umu.com}
\equalcont{These authors contributed equally to this work.}

\author[2]{\fnm{Andrian} \sur{Uihlein}}\email{andrian.uihlein@fau.de}
\equalcont{These authors contributed equally to this work.}

\author[2,3]{\fnm{Lukas} \sur{Pflug}}\email{lukas.pflug@fau.de}

\author[1,4]{\fnm{Eddie} \sur{Wadbro}}\email{eddie.wadbro@kau.se}

\affil*[1]{\orgdiv{Department of Computing Science}, \orgname{Ume{\aa} University}, \orgaddress{\city{Ume{\aa}}, \postcode{SE--90187}, \country{Sweden}}}

\affil[2]{\orgdiv{Department of Mathematics, Chair of Applied Mathematics (Continuous Optimization)}, \orgname{Friedrich-Alexander-Universität Erlangen-Nürnberg (FAU)}, \orgaddress{\street{Cauerstraße 11}, \city{Erlangen}, \postcode{91058}, \state{State}, \country{Germany}}}

\affil[3]{\orgdiv{Competence Unit for Scientific Computing (CSC)}, \orgname{Friedrich-Alexander-Universität Erlangen-Nürnberg (FAU)}, \orgaddress{\street{Martensstrasse 5a}, \city{Erlangen}, \postcode{91058}, \country{Germany}}}

\affil[4]{\orgdiv{Department of Mathematics and Computer Science}, \orgname{Karlstad University}, \orgaddress{\city{Karlstad}, \postcode{SE--65188}, \country{Sweden}}}

\abstract{This paper focuses on the topology optimization of a broadband acoustic transition section that connects two cylindrical waveguides with different radii. The primary objective is to design a transition section such that it maximizes the transmission of a planar acoustic wave while ensuring the planarity of the transmitted wave. 
Helmholtz equation is used to model linear wave propagation in the device. We utilize the finite element method to solve the state equation on a structured mesh of square elements. Subsequently, a material distribution topology optimization problem is formulated to optimize the distribution of sound-hard material in the transition section.
We employ two different gradient-based approaches to solve the optimization problem: namely, a deterministic approach using the method of moving asymptotes (MMA), and a stochastic approach utilizing both stochastic gradient (SG) and continuous stochastic gradient (CSG) methods.
A comparative analysis is provided among these methodologies concerning the design feasibility and the transmission performance of the optimized designs, and the computational efficiency.
The outcomes highlight the effectiveness of stochastic techniques in achieving enhanced broadband acoustic performance with reduced computational demands and improved design practicality. The insights from this investigation demonstrate the potential of stochastic approaches in acoustic applications, especially when broadband acoustic performance is desired.}




\maketitle

\section{Introduction}\label{sec1}
In many acoustic applications, acoustic waveguides connect various parts of a device from the transducer to the radiating aperture or receiver~\cite{Rutsch2022,Haugwitz2022}. This includes the transmission of acoustic waves between parts with different geometry, dimensions, and acoustic impedances. In most cases, a minimum reflection and alteration of the incoming signal is desired, which can be achieved by using an impedance matching section~\cite{Wadbro2014,Robertson2019}.

One of the earliest works on acoustic transition sections is the study by Kirby~\cite{Kirby2008}, who used numerical simulations to model acoustic wave propagation in two waveguides connected through a transition section.
Wadbro~\cite{Wadbro2014} used a material distribution topology optimization method to design a transition section between two cylindrical waveguides with different radii to achieve impedance matching. 
Robertson et al.~\cite{Robertson2019} considered a similar problem of impedance matching between two cylindrical waveguides, although they used a two-neck Helmholtz resonator as the transition section. They showed that perfect impedance matching can be achieved by tuning the dimensions of the Helmholtz resonator. However, this approach has bandwidth limitations with impedance matching achieved only at a narrow range of frequencies close to the resonance frequency of the resonator.
Cao et al.~\cite {Cao2020} took a different approach by considering a two-dimensional acoustic transformation section with an impedance-tunable transformed medium. They showed that desirable broadband impedance matching can be achieved in this way, though in practice, it is very difficult to set a transformation medium with acoustic properties changing according to a given function~\cite{Cao2020}. 

Here, we use the material distribution topology optimization method, also known as density-based topology optimization, to design an acoustic transition section for impedance matching. This is a common method in computational design optimization for acoustic problems~\cite{Wadbro2006,Duehring2008,Bokhari2021,Yoon2020}. To obtain a final topology that is favorable for a broad range of frequencies, there are two fundamentally different approaches:

First, the problem can be viewed as a topology optimization task for infinitely many load cases. Using an appropriate discretization scheme, the problem can be transformed into a deterministic multi-load formulation, where the objective function admits a finite sum structure~\cite{BendsoeMultiload,Multiload2,Multiload3}, for which various deterministic optimization schemes have been explored. In this contribution, the resulting optimization problems are solved via the method of moving asymptotes (MMA)~\cite{Svan87}. We investigate how the solution depends on the  number of frequencies considered in the optimization.

Changing perspective, the original objective can also be considered as a robust optimization problem~\cite{Robust1,Robust2,Robust3}, where the broadband frequency range models an underlying uncertainty of load cases. While many approaches for robust optimization again rely on discretizations or series expansions of the full objective, we specifically choose two methods from stochastic optimization which do not follow this philosophy. Namely, we optimize the transition section using the stochastic gradient descent~\cite{Monro1951} and the continuous stochastic gradient descent~\cite{CSG1,CSG2,CSG3} method. Both of these approaches represent probabilistic, sample-based optimization schemes.

The remainder of this paper is structured as follows. In Section~\ref{section:PS}, we introduce the problem setup, including the geometric configuration and the governing equations. We also discuss the discretization of the state equations using the finite element method.
Next, we introduce the objective function for the design optimization problem aiming for a broadband transition.
For the different optimization approaches considered in this work, we present the corresponding results of our numerical experiments in Section~\ref{section:NE}.
Finally, in Section~\ref{section:CD}, we present a concluding discussion that includes a comprehensive comparison of the results obtained by the different approaches. We summarize the main findings and provide insights into the strengths and limitations of each method.

\section{Problem Statement}\label{section:PS}
\begin{figure}
        \centering
	\includegraphics[width=0.9\linewidth]{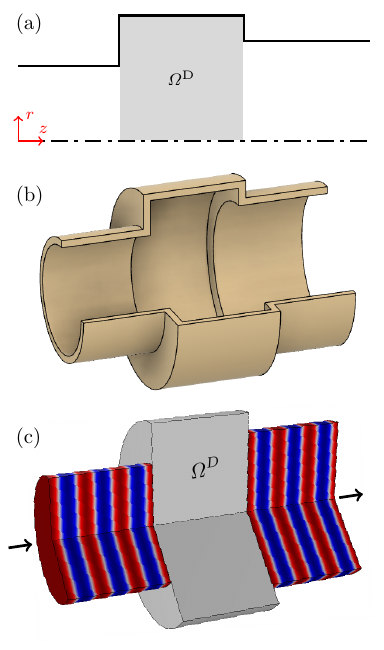}
	\caption{(a)~Axi-symmetric setup with cylindrical design domain $\Omega^\text{D}$ in the middle and one cylindrical waveguide on both sides. (b)~A 3D visualization of the setup. (c)~The targeted wave propagation characteristics in the waveguides, where the planar incoming wave in the left pipe continues to propagate as a planar wave to the right pipe.}%
	\label{fig:problem_setup}
\end{figure}

Consider the cylindrical setup illustrated in Fig.~\ref{fig:problem_setup}, consisting of two semi-infinite pipes connected by a transition section. Assume a planar acoustic wave propagating from left to right in the left pipe. As this incoming wave propagates in the transition section $\Omega^\text{D}$, highlighted in grey, a part of the wave will propagate through the transition section to the right pipe, while another part will be reflected back to the left pipe. 
By optimizing the distribution of sound-hard material in the transition section $\Omega^\text{D}$, this study aims to ensure that the planar incoming wave in the left pipe continues to propagate as a planar wave to the right pipe, despite the change in the diameter of the waveguide. A similar problem was previously considered by Wadbro~\cite{Wadbro2014} for a narrower range of frequencies. In this study, we aim to extend the analysis to a broader range of frequencies, which is of practical importance for a wide range of acoustic applications. Additionally, we compare two distinct approaches to solving the optimization problem: the deterministic approach and the stochastic approach. Here and throughout this article, the waves that propagate away from the transition section are termed \emph{outgoing} waves, and the waves that propagate towards the transition section are termed \emph{incoming} waves. We note that the outgoing waves are further divided into \emph{reflected} waves, traveling to the left in the left pipe, and \emph{transmitted} waves, traveling to the right in the right pipe.

\subsection{Mathematical model}
\begin{figure}
	\includegraphics[width=1\linewidth]{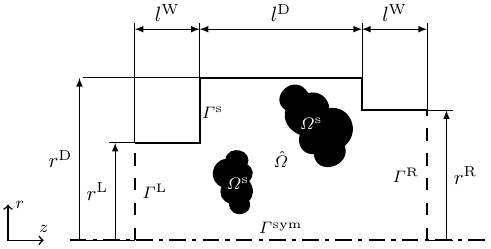}
	\caption{The computational domain $\Omega$ comprising the design domain $\Omega^\text{D}$ and two truncated waveguides.
  An arbitrary sound-hard material distribution inside $\Omega^\text{D}$ is given as $\Omega^\text{s}$\!. The remainder of the setup $\hat{\Omega}=\Omega\setminus\Omega^\text{s}$ is filled with air.
The boundaries are divided into an axi-symmetric axis $\Gamma^\text{sym}$ (dash-dotted line), the sound-hard walls $\Gamma^\text{s}$ (solid line) and the boundaries $\Gamma^\text{L}$ and $\Gamma^\text{R}$ (dashed lines) of the artificially truncated waveguides.}%
	\label{fig:computational_domain}
\end{figure}

In this study, we consider linear wave propagation in the cylindrically symmetric setup, illustrated in Fig.~\ref{fig:problem_setup}(a). Specifically, we let ${P(x,t)=Re\lc p(x)e^{i\omega t}\rc}$ denote the time-harmonic pressure in the air-filled region. Under these assumptions, the complex pressure amplitude $p(x)$ satisfies the Helmholtz equation in cylindrical coordinates. That is,
\begin{equation}\label{e:helmholtz}
-\nabla \cdot (r \nabla p) - k^2 rp = 0, \quad \text{in}\; \hat \Omega,
\end{equation}
which describes the behavior of sound waves in the system. Therein, $\hat{\Omega}$ is the air-filled region of the setup, and the wave number $k= \omega / c$ is determined by the speed of sound in air $c$ and the angular frequency $\omega$. 
As mentioned earlier, the design domain $\Omega^\text{D}$ (transition section) may partly be filled with sound-hard material. Fig.~\ref{fig:computational_domain} shows the air-filled domain $\hat{\Omega}$ and an arbitrary distribution of sound-hard material $\Omega^s$ in the design domain. By varying the distribution of sound-hard material within the design domain, we can explore how this affects the propagation of sound waves through the transition section and identify optimal designs that minimize wave reflection and ensure planar wave transmission.

To numerically approximate the infinitely long pipes on both sides of the design domain, we truncate them and use DtN (Dirichlet-to-Neumann) non-reflecting boundary conditions at the artificial boundaries $\Gamma^\text{L}$ and $\Gamma^\text{R}$, as illustrated in Fig.~\ref{fig:computational_domain}. Regarding more details about this type of artificial boundary conditions, we refer the reader to the book by Ihlenburg~\cite{Ih1998} and the appendix of the article by Wadbro~\cite{Wadbro2014}. Considering these artificial boundary conditions, we obtain the boundary-value problem
\begin{subequations}\label{eq:CstateEq}
\begin{alignat}{2}
 -\nabla \cdot(r\nabla p) - k^2 rp &= 0, &\quad &\text{in}\; \hat \Omega, \label{eq:RadHelm} \\
 \frac{\partial p}{\partial n} &= 0, &\quad &\text{on}\; \Gamma^\text{s}, \label{eq:SoundHard} \\
 \frac{\partial p}{\partial n} - \text{DtN}(p) &= 2\mathrm{i}k, &\quad &\text{on}\; \Gamma^\text{L}, \label{eq:LeftDtn}\\
 \frac{\partial p}{\partial n} - \text{DtN}(p) &= 0, &\quad &\text{on}\; \Gamma^\text{R}. \label{eq:RightDtn}
\end{alignat}
\end{subequations}
Conditions~\eqref{eq:LeftDtn} and \eqref{eq:RightDtn} ensure that all the outgoing waves are perfectly absorbed. Furthermore, condition~\eqref{eq:LeftDtn} also specifies an incoming planar wave with unit amplitude at $\Gamma^\text{L}$.
By multiplying equation~\eqref{eq:RadHelm} with a test function $q$ and integrating over the domain $\hat \Omega$, the variational form of boundary-value problem~\eqref{eq:CstateEq} can be written as follows.
\begin{multline}\label{VrEq}
\text{Find $p \in H^1 (\hat \Omega)$ such that:}\\
\int_{\hat \Omega} r \nabla q \cdot \nabla p\dd\Omega- k^2 \int_{\hat \Omega} r q p\dd\Omega \\
- \int_{\Gamma^\text{L}} rq \text{DtN}(p)\dd\Gamma - \int_{\Gamma^\text{R}} rq \text{DtN}(p)\dd\Gamma\\
 =  2\mathrm{i}k\int_{\Gamma^\text{L}} rq\dd\Gamma, \quad \forall q \in H^1 (\hat{\Omega}).
\end{multline}
For a given distribution of solid material in the design domain and a given shape of region $\Omega^s$, the solution $p$ to equation~\eqref{VrEq} shows the distribution of the complex pressure in $\hat \Omega$. Following a standard approach in topology optimization, we define a material indicator function $\alpha$ such that $\alpha \equiv 0$ in $\Omega^s$ and $\alpha \equiv 1$ in $\hat \Omega$. 
Using this function, we extend the integration domain of the domain integrals in variational formulation~\eqref{VrEq} from $\hat{\Omega}$ to ${\Omega = \hat{\Omega} \cup \Omega^s}$. We note that the computational domain $\Omega$ consists of both the air and the solid regions. The resulting reformulation of variational formulation~\eqref{VrEq} is then given by
\begin{multline}\label{ExVrEq}
\text{Find $p \in H^1 (\Omega)$ such that}\\
\int_{\Omega} \alpha r \nabla q \cdot \nabla p\dd\Omega - k^2 \int_{\Omega} \alpha r q p\dd\Omega \\
- \int_{\Gamma^\text{L}} rq \text{DtN}(p)\dd\Gamma - \int_{\Gamma^\text{R}} rq \text{DtN}(p)\dd\Gamma \\
= 2\mathrm{i}k\int_{\Gamma^\text{L}} rq\dd\Gamma, \quad \forall q \in H^1 (\Omega).
\end{multline}
The solution $p$ to the variational formulation~\eqref{ExVrEq} represents the distribution of complex pressure in the waveguide, given a design of solid scatter in $\Omega^\text{D}$, described by a material indicator function $\alpha$.

\subsection{Discretization}
We use the finite element method to discretize and numerically solve problem~\eqref{ExVrEq} on a structured grid of square elements. Let $V$ be a finite element functional space consisting of continuous and bi-quadratic functions on each element, and let $\varphi_j$, ${j=1, 2, \ldots, N}$ be bi-quadratic shape functions, where $N$ is the number of degrees of freedom. Thus, $V=\text{span}\lc\varphi_1, \varphi_2, \ldots, \varphi_N\rc$. We approximate the complex pressure $p$ and the test function $q$ by $p_h \in V$ and $q_h \in V$, respectively. Additionally, we approximate the material indicator function $\alpha$ with an element-wise constant function $\alpha_h$. Using the above definitions and approximations, we obtain the discretized version of problem~\eqref{ExVrEq} below.
\begin{multline}\label{e:var_dis}
\text{Find $p_h\in\,V$ such that}\\
  \int_\Omega \alpha_h r \nabla q_h  \cdot \nabla p_h\dd\Omega - k^2 \int_\Omega \alpha_h r q_h p_h\dd\Omega \\
  - \int_{\Gamma^\text{L}} r q_h \text{DtN}_h(p_h)\dd\Gamma
    -\int_{\Gamma^\text{R}} r q_h \text{DtN}_h(p_h)\dd\Gamma
    \\ = 2\mathrm{i}k\int_{\Gamma^\text{L}} r q_h\dd\Gamma, \quad \forall q_h\in V,
\end{multline}
where $\text{DtN}_h$ represents semi-discrete Dirichlet-to-Neumann type boundary operators on $\Gamma^\text{L}$ and $\Gamma^\text{R}$~\cite[Appendix~A]{Wadbro2014}.
The algebraic or matrix formulation of problem~\eqref{e:var_dis} reads
\begin{equation}
	\left(\mathbf{K}(\boldsymbol{\alpha}) - k^2\mathbf{M}(\boldsymbol{\alpha}) - \mathbf{B}^\text{L} -  \mathbf{B}^\text{R}\right) \mathbf{p} = 2ik \mathbf{M}^\text{L} \mathbf{\mathds{1}_N},
	\label{eq:goveq}
\end{equation}
where $\mathbf{p}=\left[p_1,p_2,\ldots,p_{N}\right]^T$ is the vector of nodal values of the complex acoustic pressure amplitude, $\boldsymbol{\alpha}=\left[\alpha_1,\alpha_2,\ldots,\alpha_{N^\text{D}}\right]^T$ is the vector that holds the element values of $\alpha_h$ (with $N^\text{D}$ denoting the number of elements in $\Omega^\text{D}$) and $\mathds{1}_N = [1,1, \dots, 1]^T$ is a vector of length $N$. Also, the $N\times N$ stiffness $\mathbf{K}$, mass $\mathbf{M}$, and boundary mass $\mathbf{M}^\text{L}$ matrices have components
\begin{subequations}
\begin{alignat}{1}
    K_{i j} &= \int_\Omega \alpha_h r \nabla\varphi_i \cdot \nabla\varphi_j\dd\Omega,\\
    M_{i j} &= \int_\Omega \alpha_h r \varphi_i \varphi_j\dd\Omega,\\
    M^\text{L}_{i j} &= \int_{\Gamma^\text{L}}\; r\varphi_i \varphi_j \dd\Gamma, \label{eq:ML}
\end{alignat}
\end{subequations}
respectively. The boundary matrices $\mathbf{B}^\text{L}$ and $\mathbf{B}^\text{R}$ represent the non-reflecting boundary conditions at $\Gamma^\text{L}$ and $\Gamma^\text{R}$, respectively. A detailed derivation of $\mathbf{B}^\text{L}$ and $\mathbf{B}^\text{R}$ is provided in a previous work by Wadbro~\cite[Appendix A]{Wadbro2014}.

\subsection{Power of outgoing waves}
Let Helmholtz equation~\eqref{e:helmholtz} govern the distribution of the complex pressure $p$ in the two semi-infinite pipes on the left and right side of the transmission section illustrated in Fig.~\ref{fig:problem_setup}(a) with sound-hard boundary condition~\eqref{eq:SoundHard} on the solid walls. 
Using the separation of variables, the general solution for $p$ in the left and right pipes reads
\begin{subequations}\label{eq:SOVP}
\begin{align}
      p^\text{L} &= \sum_m f_m(r) \bigg(A^\text{L}_m \mathrm{e}^{\mathrm{i} k_m (z^\text{L}-z)} + B^\text{L}_m \mathrm{e}^{\mathrm{i} k_m (z-z^\text{L})}\bigg),\\
      p^\text{R} &= \sum_m f_m(r) \bigg(A^\text{R}_m \mathrm{e}^{\mathrm{i} k_m (z-z^\text{R})} + B^\text{R}_m \mathrm{e}^{\mathrm{i} k_m (z^\text{R}-z)}\bigg),
\end{align}
\end{subequations}
respectively, where $z^\text{L}$ and $z^\text{R}$ are the position of $z$-axis on $\Gamma^\text{L}$ and $\Gamma^\text{R}$, the functions $f_m(r)$ are the modes at left and right waveguide, $\mathrm{e}$ is the base of the natural logarithm, $\mathrm{i}$ is the imaginary unit and $A^\text{L}_m$ and $B^\text{L}_m$ are complex constants that determine the amplitude of incoming and outgoing waves at the left waveguide, respectively.
Similarly, $A^\text{R}_m$ and $B^\text{R}_m$ are complex constants that determine the amplitude of incoming and outgoing waves at the right waveguide, respectively. Lastly, the constants $k_m$ are the so-called reduced wave numbers.
In the continuous case, we have an infinite number of modes $m=0,1,2,\ldots$, but only the modes with real-values of $k_m$ are propagating modes and the ones with imaginary $k_m$ will decay exponentially according to the equations~\eqref{eq:SOVP}. These modes are known as evanescent modes.

The mode functions $f_m(r)$ should satisfy the following one-dimensional eigenvalue problem in the radial direction on the boundaries $\Gamma^\text{L}$ and $\Gamma^\text{R}$:
\begin{equation}
  -\frac{\partial}{\partial r}\bigg(r \frac{\partial f}{\partial r}\bigg) = \lambda r f,
  \quad
  \frac{\partial f}{\partial r} \bigg\rvert_{r=0} =\frac{\partial f}{\partial r} \bigg\rvert_{r=W}=0,
  \label{eq:eigenproblem}
\end{equation}
where $W$ is the radius of the pipe. 
In the continuous case, it is well known that the functions $f_m$ are so-called Bessel functions.
For our numerical treatment, we can extrapolate the numerical solution on $\Gamma^\text{L}$ and $\Gamma^\text{R}$ to any point in the waveguides using expansion~\eqref{eq:SOVP} and viewing the problem continuously in the lengthwise direction and discretely in the radial direction.
Note that for a given finite element discretization of equation~\eqref{eq:eigenproblem}, the number of modes that are representable in the discretized case $M$ equals the number of basis functions with support on the boundary.
Let $f_m^h$ be an eigenfunction (mode) corresponding to eigenvalue $\lambda_m$, where $m=0,1,\ldots, M$.
Since the complex pressure $p$ satisfies Helmholtz equation~\eqref{e:helmholtz}, for the reduced wave number we have
\begin{equation}\label{eq:k_m}
k_{m}^2 = k^2 - \lambda_{m}.
\end{equation}
Recall that $f_m^h$ is a propagating mode if its corresponding reduced wave number is real. The smallest eigenvalue in solving problem~\eqref{eq:eigenproblem} is 0, which corresponds to the planar wave. 
So $\lambda_{0}=0$ and $k_0=k$, and thus, the planar wave mode is always a propagating mode. Moreover, for a given frequency $f$, there is a finite number $M^p$ of propagating modes satisfying the condition $\lambda_m \leq k^2$. The number of propagating modes depends on the frequency of the wave and the radius of the pipe. Thus, for the different radii of the left and right pipes, we may have a different number of propagating modes at $\Gamma^\text{L}$ and $\Gamma^\text{R}$, which we denote by $M_L^p$ and $M_R^p$, respectively.

In the discretized case, the solution at the boundaries $\Gamma^\text{L}$ and $\Gamma^\text{R}$, where $z^\text{L}-z=0$ and $z-z^\text{R}=0$, respectively, reads 
\begin{subequations}\label{eq:BSOVP}
\begin{align}
      p\bigg\rvert_{\Gamma^\text{L}} &= \sum_0^{M_L^p} f_m^h(r) \bigg(A^\text{L}_m + B^\text{L}_m \bigg),\\
      p\bigg\rvert_{\Gamma^\text{R}} &= \sum_0^{M_R^p} f_m^h(r) \bigg(A^\text{R}_m + B^\text{R}_m \bigg).
\end{align}
\end{subequations}

From now onward, we occasionally use the superscript $X$ in an expression to represent either $\text{L}$ for the left waveguide or $\text{R}$ for the right waveguide. Note that the corresponding statement holds in both cases of replacing $X$ by either $\text{L}$ or $\text{R}$ referring to the left and right waveguides, respectively. Let $f_n^h$ be the $n$th propagating mode. Then, we have
\begin{equation}\label{eq:ortho_B_m}
\begin{split}
    \int_{\Gamma^X} r f_n^h p &= \sum_{0}^{M^p} \int_{\Gamma^X} r f_n^h \bigg( \big(A_m^X + B_m^X \big) f_m^h \bigg) 
	\\ &= A_n^X + B_n^X,
\end{split}
\end{equation}
where the first equality follows from substituting $p\big\rvert_{\Gamma^\text{X}}$ from equation~\eqref{eq:BSOVP} into the first expression and the second equality follows from the orthonormality of modes. Considering $\mathbf{v}_n^X$ to be the $N \times 1$ vector representing the nodal values of the discrete mode $f_n^h$ on the boundary nodes at $\Gamma^X$ (note that all other entries of $\mathbf{v}_n^X$ corresponding to the internal nodes are zero), equation~\eqref{eq:ortho_B_m} in matrix form reads
\begin{equation}\label{eq:Amplitudes}
	 \lp\mathbf{v}_n^X\rp^T\mathbf{M}^\text{X} \mathbf{p}= A_n^X + B_n^X,
\end{equation}
where $\mathbf{M}^\text{X}$ is the boundary mass matrix as defined in equation~\eqref{eq:ML} at $\Gamma^X$, and $\mathbf{p}$ is the nodal values of the complex pressure. Thus, for a given solution $\mathbf{p}$ to problem~\eqref{eq:goveq}, we can recover the complex amplitudes of the incoming and outgoing waves for each propagating mode at $\Gamma^X$, using equation~\eqref{eq:Amplitudes}.

In this study, we only consider the case where we have a planar incoming wave with unit amplitude at $\Gamma^\text{L}$. Therefore, we can rewrite equation~\eqref{eq:Amplitudes} as
\begin{equation}\label{eq:Amplitudes_fop}
\begin{split}
    \lp\mathbf{v}_0^\text{L}\rp^T\mathbf{M}^\text{L} \mathbf{p} &= 1 + B_0^\text{L},\\
    \lp\mathbf{v}_m^\text{L}\rp^T\mathbf{M}^\text{L} \mathbf{p} &= B_m^\text{L}, \quad m=1,2,\ldots,M_L^p,\\
    \lp\mathbf{v}_n^\text{R}\rp^T\mathbf{M}^\text{R} \mathbf{p} &= B_n^\text{R}, \quad n=0,1,\ldots,M_R^p.\\
\end{split}
\end{equation}
Recall that $M_L^p$ and $M_R^p$ are the highest propagating modes in the left and right pipes, respectively.

The power of a propagating wave is proportional to the square of its amplitude and its corresponding reduced wave number. Defining the normalized power of outgoing waves as the power of outgoing wave divided by the power of the unit-amplitude incoming wave and considering equation~\eqref{eq:Amplitudes_fop} to compute the amplitude of outgoing waves, we have
\begin{equation}\label{eq:Powers_fopL}
\begin{split}
    P_m^\text{L}=
        \begin{cases}
            \bigg \lvert\lp\mathbf{v}_0^\text{L}\rp^T\mathbf{M}^\text{L} \mathbf{p}-1 \bigg\rvert^2 \; &\text{if} \; m=0,\\[10pt]
            \frac{k_m}{k}\bigg\lvert\lp\mathbf{v}_m^\text{L}\rp^T\mathbf{M}^\text{L} \mathbf{v}\bigg\rvert^2 \; &\text{if} \; m=1,\ldots,M_L^p,\\
	\end{cases}
 \end{split}
\end{equation}
and
\begin{equation}\label{eq:Powers_fopR}
    P_n^\text{R}=
        \frac{k_n}{k}\bigg\lvert\lp\mathbf{v}_n^\text{R}\rp^T\mathbf{M}^\text{R} \mathbf{p}\bigg\rvert^2 \; \text{for} \; n=0,1,\ldots,M_R^p.
\end{equation}
Here $P_m^\text{L}$ and $P_n^\text{R}$ are the normalized power of the outgoing waves of mode $m$ and $n$ at $\Gamma^\text{L}$ and $\Gamma^\text{R}$, respectively. A detailed derivation of the power of outgoing waves for a given amplitude in the discretized case is provided by Wadbro~\cite[Appendix B]{Wadbro2014}. Here and throughout this article, whenever the term \emph{power of outgoing wave} is used, it means the \emph{normalized} power of the outgoing wave by dividing the power of the outgoing wave with the power of a unit-amplitude incoming wave imposed at $\Gamma^\text{L}$.

\subsection{Objective function}
As mentioned earlier, the aim of this study is to design the transmission section in Fig.~\ref{fig:problem_setup} to (i) maximize the transmission and (ii) ensure that the transmitted wave is planar as illustrated in Fig.~\ref{fig:problem_setup}(c). To achieve this, we minimize the sum of the power of all outgoing waves except for the planar wave to the right over the targeted range of frequencies $\mathscr{F}$. Thus, the primary objective function can be written as
\begin{equation}\label{eq:Objective}
    J_p(\boldsymbol{\alpha}) = \frac{1}{\lvert \mathscr{F} \rvert}\int_{\mathscr{F}} \lp \sum^{M^{p}_L\lp f \rp}_{m=0}P_m^\text{L} + \sum^{M^{p}_R\lp f \rp}_{n=1}P_n^\text{R}\rp \dd \mathscr{F}.
\end{equation}
Note that we have normalized the objective by the length of the targeted frequency range $\lvert \mathscr{F} \rvert$. Note that for each frequency, we need to calculate $M^{p}_L\lp f \rp$ and $M^{p}_R\lp f \rp$, the number of propagating modes at $\Gamma^\text{L}$ and $\Gamma^\text{R}$, respectively.

For binary values of $\alpha_h=\lc0,1\rc$, the optimization problem with objective function~\eqref{eq:Objective} is a large-scale non-linear integer optimization problem. Also, if $\alpha_h=0$ in some of the elements, then the system matrix in equation~\eqref{eq:goveq} becomes singular.
To solve the numerical and mathematical issues that arise when solving this problem, a standard approach in topology optimization is to relax the binary value constraint and let $\alpha_h$ take values in the range $[\epsilon,1]$, where $\epsilon$ is a small number~\cite{Wadbro2006,DuJeSi2008,Wadbro2014,KaWaBe2015,Bokhari2021}. Moreover, we aim for a pure solid ($\alpha_h=\epsilon$) or air ($\alpha_h=1$) final design. Thus, we use a combination of filtering and penalty methods to suppress the intermediate values. The non-linear density filters used in the numerical experiments also ensure a size control on the solid region in the design~\cite{HaWaHaBe2018,Sigmund2007,HaWa17,Hagg2018,Bokhari2021}. Let ${\mathbf{d}=\left[d_1,d_2,\ldots,d_{N^\text{D}}\right]^T}$ be the vector of design variables before filtering. 
Thus, we define the $N^\text{D} \times 1$ vector ${\boldsymbol{\alpha} :=\mathcal{F}\lp \mathbf{d} \rp}$, where $\mathcal{F}$ is a filter operator.
To further suppress the intermediate values of the design variables, we add a standard quadratic penalty term~\cite{Allaire1993,Borrvall2001,Wadbro2014,Bokhari2021} to the primary objective function~\eqref{eq:Objective} and thus, we define the objective function for the numerical experiments as follows:
\begin{equation}\label{eq:FinalObjective}
\begin{split}
     J(\boldsymbol{\alpha})&=J_p(\boldsymbol{\alpha}) + \frac{\gamma }{\lvert \Omega^\text{D} \rvert}\int_{\Omega^\text{D}} (\alpha_h-\epsilon)(1-\alpha_h) \\
     &= \frac{1}{\lvert \mathscr{F} \rvert} \int_{\mathscr{F}} \lp \sum^{M^{p}_L\lp f \rp}_{m=0}P_m^\text{L} + \sum^{M^{p}_R\lp f \rp}_{n=1}P_n^\text{R} \rp \dd \mathscr{F} \\
     &+ \frac{\gamma }{N^\text{D}} \sum^{N^\text{D}}_{k=1}(\alpha_k-\epsilon)(1-\alpha_k),
\end{split}
\end{equation}
where $\gamma$ is the penalty parameter, $\lvert \Omega^\text{D} \rvert$ denotes the size of the design domain.

\section{Numerical experiments}\label{section:NE}
In our numerical experiments, we consider the setup illustrated in Fig.~\ref{fig:computational_domain} with the following dimensions: The radius and length of the design region $\Omega^\text{D}$ is ${r^\text{D}=50\,\text{mm}}$ and ${l^\text{D}=50\,\text{mm}}$, respectively. The radius and length of the truncated waveguides are ${r^\text{L}=30\,\text{mm}}$, ${r^\text{R}=40\,\text{mm}}$, and ${l^\text{W}=20\,\text{mm}}$, respectively. We aim to maximize the transmission of the planar incoming wave in the frequency range of 4--16\,kHz, ensuring that the transmitted wave is also planar.
We discretize the computational domain into a structured grid of square elements with a uniform mesh size of $h=0.25\;\text{mm}$, resulting in 250,721 degrees of freedom for the finite element discretization.
To solve the optimization problem, we employ three different optimization algorithms: the MMA method~\cite{Svan87}, the stochastic gradient (SG) method~\cite{Monro1951}, and the continuous stochastic gradient (CSG) method~\cite{CSG1,CSG2,CSG3}. 

We define the performance of a given design at frequency $f$ as the normalized power of the outgoing planar wave, computed using expression~\eqref{eq:Powers_fopR}, as follows:
\begin{equation}\label{eq:Performance}
    \text{Performance}=P_0^\text{R}(f)=
    \bigg\lvert\lp\mathbf{v}_0^\text{R}\rp^T\mathbf{M}^\text{R} \mathbf{p}\bigg\rvert^2.
\end{equation}
To evaluate the performance of the optimized designs, we use a boundary-fitted mesh for the final designs in the Acoustics Modules in COMSOL Multiphysics. 
The performance of different designs are compared over a range of frequencies from 4\,kHz to 16\,kHz with a step size of 20\,Hz.

\subsection{MMA approach}
To solve the optimization problem considering objective function~\eqref{eq:FinalObjective} using the MMA method, we approximate the integral over the range of targeted frequencies using the function values of the integrand at just a few frequencies. 
Thus, we discretize the optimization problem as
\begin{multline}\label{eq:MMAOptProb}
    \min_{\mathbf{d}\in \mathcal{A}} \quad \sum^{Q}_{i=1} \lp \sum^{M^{p}_L\lp f_i \rp}_{m=0}P_m^\text{L} + \sum^{M^{p}_R\lp f_i \rp}_{n=1}P_n^\text{R} \rp \\
    + \frac{\gamma}{N^\text{D}} (\mathcal{F}(\mathbf{d})-\epsilon\mathbf{\mathds{1}_{N^\text{D}}})^T(\mathbf{\mathds{1}_{N^\text{D}}}-\mathcal{F}(\mathbf{d})), 
\end{multline}
where $Q$ is the number of frequencies subject to optimization, $\mathbf{\mathds{1}_{N^\text{D}}}$ is $N^\text{D} \times 1$ vector with all entries equal to 1, and $\mathcal{A} =  \{\mathbf{d}\in \mathbb{R}^{N^\text{D}} \mid \epsilon \leq d_i \leq 1 \;\forall\; i  \}$ is the set of admissible designs.
The scaling constant $Q^{-1}$ is neglected in expression~\eqref{eq:MMAOptProb}. This is done because the scaling between the primary objective function and the quadratic penalty term can be also tuned using the penalty parameter $\gamma$.
Note that by increasing the number of frequencies subject to optimization $Q$, we can improve the approximation used to discretize objective function~\eqref{eq:FinalObjective}.
To solve optimization problem~\eqref{eq:MMAOptProb}, we utilize the least squares formulation of the MMA approach described by Svangberg~\cite{Svan87,Svan2002}. 
Thus, we need the sensitivity information for each part of the objective function in optimization problem~\eqref{eq:MMAOptProb}. 
The computation of sensitivities for the quadratic penalty term can be readily performed for a given filter $\mathcal{F}$. However, the task of computing the gradient of the power of outgoing modes with respect to the design variables poses a challenge. This process involves determining the gradient of the amplitudes of the outgoing modes with respect to the design variables, as per equations~\eqref{eq:Powers_fopL} and~\eqref{eq:Amplitudes_fop}. Notably, the power of a propagating mode is proportional to the square of its amplitude.
The sensitivity computations are done analytically using the adjoint variable method, which is a powerful technique for computing the gradient of a function that depends on the solution of a partial differential equation with respect to the design variables. A detailed derivation of the design sensitivities is provided in the appendix.

The design problem has 40,000 design variables, that is, the number of elements in the design domain $\Omega^\text{D}$.
We set the lower bound for the design variable as $\epsilon=10^{-8}$, the filter radius to $1\,\text{mm}$, and use a so-called continuation approach for the penalty parameter.
That is, we solve problem~\eqref{eq:MMAOptProb} for a sequence of increasing penalty parameters $\gamma_i = 10^{i}$, $i=0,1,\ldots,5$, using the previously computed solution as the initial design. 
The aim is to gradually move the optimizer's focus from the acoustic performance of the device towards obtaining a black-and-white final design.
This approach ensures an optimized layout with sharp solid--fluid boundaries, free of any intermediate values of the material indicator functions~\cite{Wadbro2006}.
To solve optimization problem~\eqref{eq:MMAOptProb}, we need to consider a sufficiently large number of frequencies $Q$ to get broadband acoustic performance. However, increasing the number of frequencies subject to optimization will increase the number of times we need to solve the state equation~\eqref{eq:goveq}. Note that the finite element solver is the primary contributor to the computational costs in the optimizer. Consequently, increasing the number of frequencies will result in a significant increase in computational cost. Here, we consider three cases for the number of frequencies subject to optimization: 
\begin{enumerate}[leftmargin=0pt, label=\textbf{Case \Roman*}, ref=\textbf{Case~\Roman*}, itemindent=*,align=left]
\item Four equidistant frequencies in the targeted range; that is, 4\;kHz, 8\;kHz, 12\;kHz, and 16\;kHz.\label{item:Case1}
\item Seven equidistant frequencies in the targeted range; that is, 4\;kHz, 6\;kHz, \ldots, 16\;kHz.\label{item:Case2}
\item Thirteen equidistant frequencies in the targeted range; that is, 4\;kHz, 5\;kHz, \ldots, 16\;kHz.\label{item:Case3}
\end{enumerate}

\begin{figure}
        \centering
	\includegraphics[width=\linewidth]{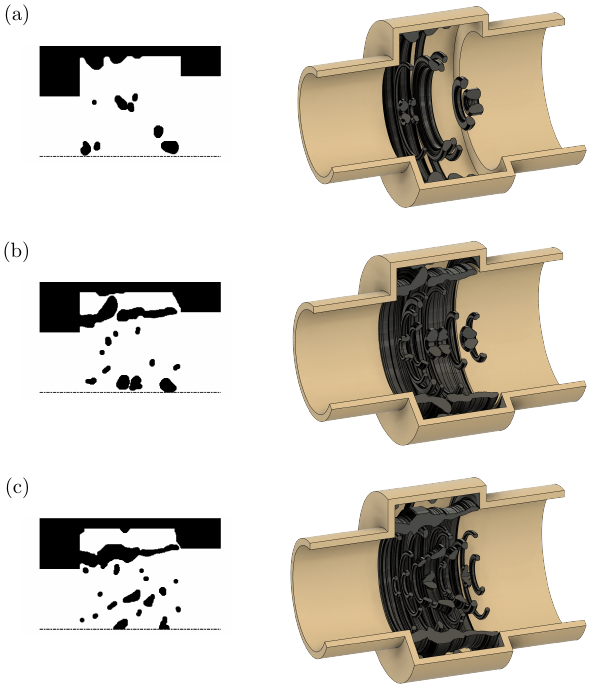}
	\caption{On the left is an axisymmetric and on the right is a 3D visualization of the optimized designs in (a) \ref{item:Case1}, (b) \ref{item:Case2}, and (c) \ref{item:Case3}}%
	\label{fig:MMADesigns}
\end{figure}

\begin{figure}
\begin{center}
    \input{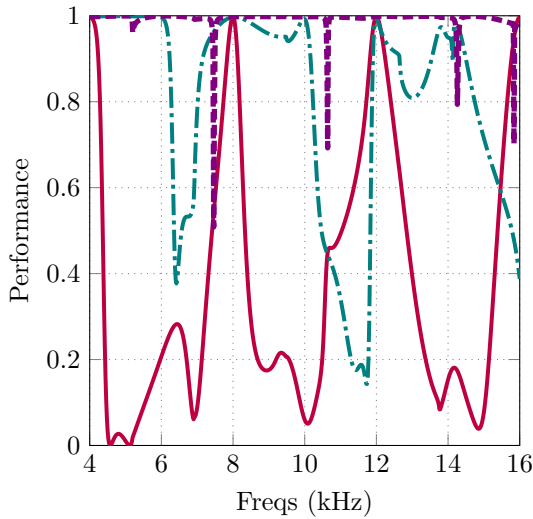}
	\caption{Performance of the optimized designs in (a) \ref{item:Case1}, the red line (b) \ref{item:Case2}, the dash-dotted teal line and (c) \ref{item:Case3}, the dashed violet line. }%
	\label{fig:MMA}
\end{center}
\end{figure}

\begin{figure}
    \textbf{Case I}\par\medskip
    \begin{center}
%
%
\definecolor{mycolor1}{rgb}{0.00000,0.44700,0.74100}%
\definecolor{mycolor2}{rgb}{0.85000,0.32500,0.09800}%
\definecolor{mycolor3}{rgb}{0.00000,0.7410000,0.32500}%
\begin{tikzpicture}

\begin{axis}[%
width=.75\linewidth,
height=.5\linewidth,
at={(0.758in,0.481in)},
scale only axis,
xmode=log,
xmin=4,
xmax=5000,
xminorticks=true,
ymode=log,
ymin=1e-3,
ymax=1,
axis background/.style={fill=white},
xmajorgrids,
xminorgrids,
ymajorgrids,
xlabel={Evaluations},
legend pos=north west
]
\addplot [color=black, line width=1.5pt, forget plot, line join=round]
  table[row sep=crcr]{%
4 0.8568\\
8 0.7723\\
12 0.6842\\
16 0.5901\\
20 0.4895\\
24 0.3839\\
28 0.2804\\
32 0.1855\\
36 0.1117\\
40 0.0645\\
44 0.036\\
48 0.019\\
52 0.0111\\
56 0.0067\\
60 0.005\\
64 0.004\\
68 0.0035\\
72 0.0033\\
76 0.0032\\
80 0.0029\\
84 0.0025\\
88 0.0023\\
92 0.0021\\
96 0.0019\\
100 0.0017\\
104 0.0016\\
108 0.0016\\
112 0.0015\\
116 0.0015\\
120 0.0015\\
124 0.0014\\
128 0.0014\\
132 0.0073\\
136 0.003\\
140 0.0022\\
144 0.0022\\
148 0.0013\\
152 0.0012\\
156 0.0011\\
160 0.0011\\
164 0.0011\\
168 0.0011\\
172 0.0011\\
176 0.0011\\
180 0.0011\\
184 0.0011\\
188 0.0011\\
192 0.0011\\
196 0.0011\\
200 0.0011\\
204 0.0011\\
208 0.0011\\
212 0.0016\\
216 0.0014\\
220 0.0014\\
224 0.0015\\
228 0.0015\\
232 0.0015\\
236 0.0016\\
240 0.0016\\
244 0.0016\\
248 0.0016\\
252 0.0016\\
256 0.0016\\
260 0.0016\\
264 0.0016\\
268 0.0016\\
272 0.0016\\
276 0.0016\\
280 0.0016\\
284 0.0016\\
288 0.0016\\
292 0.0016\\
296 0.0016\\
300 0.0016\\
304 0.0016\\
308 0.0016\\
312 0.0016\\
316 0.0016\\
320 0.0016\\
324 0.0016\\
328 0.0016\\
332 0.0016\\
336 0.0016\\
340 0.0016\\
344 0.0016\\
348 0.0016\\
352 0.0016\\
356 0.0016\\
360 0.0016\\
364 0.0045\\
368 0.003\\
372 0.003\\
376 0.0034\\
380 0.0038\\
384 0.0042\\
388 0.0045\\
392 0.0045\\
396 0.0045\\
400 0.0046\\
404 0.0046\\
408 0.0046\\
412 0.0046\\
416 0.0046\\
420 0.0046\\
424 0.0046\\
428 0.0046\\
432 0.0046\\
436 0.0046\\
440 0.0046\\
444 0.0046\\
448 0.0046\\
452 0.0046\\
456 0.0046\\
460 0.0046\\
464 0.0046\\
468 0.0046\\
472 0.0046\\
476 0.0046\\
480 0.0046\\
484 0.0046\\
488 0.0046\\
492 0.0046\\
496 0.0046\\
};

\addplot [color=purple, line width=1.5pt, forget plot, line join=round, dash pattern=on 8pt off 2pt on 1.5pt off 2pt]
  table[row sep=crcr]{%
1 0.547\\
4 0.5447\\
8 0.5404\\
12 0.5364\\
16 0.5331\\
20 0.5309\\
24 0.5291\\
28 0.5272\\
32 0.5251\\
36 0.5289\\
40 0.5383\\
44 0.5531\\
48 0.5698\\
52 0.5834\\
56 0.5936\\
60 0.6027\\
64 0.6107\\
68 0.6169\\
72 0.6215\\
76 0.6244\\
80 0.6257\\
84 0.626\\
88 0.6259\\
92 0.6256\\
96 0.625\\
100 0.6243\\
104 0.6242\\
108 0.6243\\
112 0.6245\\
116 0.625\\
120 0.6255\\
124 0.6258\\
128 0.627\\
132 0.6292\\
136 0.628\\
140 0.6262\\
144 0.6259\\
148 0.626\\
152 0.6258\\
156 0.6255\\
160 0.6262\\
164 0.626\\
168 0.6272\\
172 0.6263\\
176 0.6264\\
180 0.6266\\
184 0.6268\\
188 0.6273\\
192 0.6275\\
196 0.6275\\
200 0.6277\\
204 0.6281\\
208 0.6287\\
212 0.6282\\
216 0.6301\\
220 0.6285\\
224 0.6295\\
228 0.6286\\
232 0.6286\\
236 0.6297\\
240 0.6291\\
244 0.6282\\
248 0.6282\\
252 0.6283\\
256 0.6292\\
260 0.6291\\
264 0.628\\
268 0.6278\\
272 0.6278\\
276 0.6277\\
280 0.6275\\
284 0.6276\\
288 0.6277\\
292 0.6278\\
296 0.6279\\
300 0.628\\
304 0.6283\\
308 0.6286\\
312 0.6291\\
316 0.6295\\
320 0.6292\\
324 0.6288\\
328 0.6284\\
332 0.6282\\
336 0.6281\\
340 0.628\\
344 0.628\\
348 0.6279\\
352 0.6279\\
356 0.6279\\
360 0.6291\\
364 0.629\\
368 0.6293\\
372 0.6297\\
376 0.6301\\
380 0.6304\\
384 0.6306\\
388 0.6306\\
392 0.6307\\
396 0.6307\\
400 0.6307\\
404 0.6307\\
408 0.6307\\
412 0.6307\\
416 0.6307\\
420 0.6307\\
424 0.6307\\
428 0.6307\\
432 0.6307\\
436 0.6307\\
440 0.6307\\
444 0.6307\\
448 0.6307\\
452 0.6307\\
456 0.6307\\
460 0.6307\\
464 0.6307\\
468 0.6307\\
472 0.6307\\
476 0.6307\\
480 0.6307\\
484 0.6307\\
488 0.6307\\
492 0.6307\\
496 0.6307\\
};

\end{axis}

\end{tikzpicture}%
    \end{center}
    \textbf{Case II}\par\medskip
    \begin{center}
%
%
\definecolor{mycolor1}{rgb}{0.00000,0.44700,0.74100}%
\definecolor{mycolor2}{rgb}{0.85000,0.32500,0.09800}%
\definecolor{mycolor3}{rgb}{0.00000,0.7410000,0.32500}%
\begin{tikzpicture}

\begin{axis}[%
width=.75\linewidth,
height=.5\linewidth,
at={(0.758in,0.481in)},
scale only axis,
xmode=log,
xmin=4,
xmax=5000,
xminorticks=true,
ymode=log,
ymin=1e-3,
ymax=1,
axis background/.style={fill=white},
xmajorgrids,
xminorgrids,
ymajorgrids,
xlabel={Evaluations},
legend pos=south west
]
\addplot [color=black, line width=1.5pt, forget plot, line join=round]
  table[row sep=crcr]{%
7 0.97\\
14 0.883\\
21 0.799\\
28 0.751\\
35 0.715\\
42 0.68\\
49 0.642\\
56 0.601\\
63 0.552\\
70 0.496\\
77 0.435\\
84 0.371\\
91 0.313\\
98 0.25\\
105 0.185\\
112 0.132\\
119 0.097\\
126 0.072\\
133 0.056\\
140 0.046\\
147 0.039\\
154 0.033\\
161 0.027\\
168 0.023\\
175 0.02\\
182 0.017\\
189 0.015\\
196 0.013\\
203 0.014\\
210 0.011\\
217 0.009\\
224 0.007\\
231 0.005\\
238 0.004\\
245 0.004\\
252 0.003\\
259 0.003\\
266 0.002\\
273 0.002\\
280 0.002\\
287 0.002\\
294 0.002\\
301 0.003\\
308 0.005\\
315 0.004\\
322 0.009\\
329 0.004\\
336 0.012\\
343 0.011\\
350 0.021\\
357 0.005\\
364 0.011\\
371 0.003\\
378 0.024\\
385 0.007\\
392 0.008\\
399 0.017\\
406 0.01\\
413 0.009\\
420 0.005\\
427 0.002\\
434 0.053\\
441 0.031\\
448 0.035\\
455 0.022\\
462 0.007\\
469 0.005\\
476 0.002\\
483 0.002\\
490 0.002\\
497 0.002\\
504 0.002\\
511 0.002\\
518 0.002\\
525 0.002\\
532 0.002\\
539 0.002\\
546 0.002\\
553 0.002\\
560 0.002\\
567 0.002\\
574 0.002\\
581 0.002\\
588 0.002\\
595 0.002\\
602 0.002\\
609 0.002\\
616 0.002\\
623 0.002\\
630 0.002\\
637 0.002\\
644 0.002\\
651 0.002\\
658 0.002\\
665 0.002\\
672 0.002\\
679 0.002\\
686 0.002\\
693 0.002\\
700 0.002\\
707 0.002\\
714 0.003\\
721 0.003\\
728 0.003\\
735 0.004\\
742 0.004\\
749 0.01\\
756 0.005\\
763 0.004\\
770 0.004\\
777 0.004\\
784 0.004\\
791 0.004\\
798 0.004\\
805 0.004\\
812 0.004\\
819 0.004\\
826 0.004\\
833 0.004\\
840 0.004\\
847 0.004\\
854 0.004\\
861 0.004\\
868 0.004\\
875 0.004\\
882 0.004\\
889 0.004\\
896 0.004\\
903 0.004\\
910 0.004\\
917 0.004\\
924 0.004\\
931 0.004\\
938 0.004\\
945 0.004\\
952 0.004\\
959 0.195\\
966 0.197\\
973 0.286\\
980 0.256\\
987 0.167\\
994 0.212\\
1001 0.243\\
1008 0.239\\
1015 0.228\\
1022 0.214\\
1029 0.197\\
1036 0.177\\
1043 0.158\\
1050 0.139\\
1057 0.13\\
1064 0.125\\
1071 0.115\\
1078 0.109\\
1085 0.102\\
1092 0.096\\
1099 0.092\\
1106 0.088\\
1113 0.085\\
1120 0.082\\
1127 0.08\\
1134 0.077\\
1141 0.073\\
1148 0.07\\
1155 0.066\\
1162 0.063\\
1169 0.058\\
1176 0.054\\
1183 0.049\\
1190 0.045\\
1197 0.041\\
1204 0.037\\
1211 0.035\\
1218 0.033\\
1225 0.027\\
1232 0.024\\
1239 0.022\\
1246 0.021\\
1253 0.02\\
1260 0.019\\
1267 0.018\\
1274 0.018\\
1281 0.017\\
1288 0.017\\
1295 0.016\\
1302 0.016\\
1309 0.016\\
1316 0.016\\
1323 0.015\\
1330 0.015\\
1337 0.015\\
1344 0.016\\
1351 0.018\\
1358 0.015\\
1365 0.015\\
1372 0.015\\
1379 0.015\\
1386 0.015\\
1393 0.014\\
1400 0.014\\
1407 0.014\\
1414 0.017\\
1421 0.015\\
1428 0.015\\
1435 0.015\\
1442 0.016\\
1449 0.016\\
1456 0.016\\
1463 0.016\\
1470 0.017\\
1477 0.017\\
1484 0.018\\
1491 0.018\\
1498 0.018\\
1505 0.018\\
1512 0.019\\
1519 0.019\\
1526 0.019\\
1533 0.019\\
1540 0.019\\
1547 0.019\\
1554 0.02\\
1561 0.02\\
1568 0.02\\
1575 0.021\\
1582 0.021\\
1589 0.022\\
1596 0.022\\
1603 0.023\\
1610 0.025\\
1617 0.027\\
1624 0.028\\
1631 0.028\\
1638 0.029\\
1645 0.03\\
1652 0.031\\
1659 0.031\\
1666 0.031\\
1673 0.031\\
1680 0.031\\
1687 0.033\\
1694 0.03\\
1701 0.03\\
1708 0.03\\
1715 0.03\\
1722 0.03\\
1729 0.031\\
1736 0.031\\
1743 0.032\\
1750 0.032\\
1757 0.031\\
1764 0.032\\
1771 0.03\\
1778 0.029\\
1785 0.03\\
1792 0.057\\
1799 0.085\\
1806 0.134\\
1813 0.147\\
1820 0.208\\
1827 0.182\\
1834 0.183\\
1841 0.183\\
1848 0.183\\
1855 0.183\\
1862 0.183\\
1869 0.183\\
1876 0.183\\
1883 0.183\\
1890 0.183\\
1897 0.183\\
1904 0.183\\
1911 0.183\\
1918 0.183\\
1925 0.183\\
1932 0.183\\
1939 0.183\\
1946 0.183\\
1953 0.183\\
1960 0.183\\
1967 0.183\\
1974 0.183\\
1981 0.183\\
1988 0.183\\
1995 0.183\\
2002 0.183\\
2009 0.182\\
2016 0.182\\
2023 0.182\\
2030 0.182\\
2037 0.182\\
2044 0.182\\
2051 0.182\\
2058 0.182\\
2065 0.182\\
2072 0.182\\
2079 0.182\\
2086 0.182\\
2093 0.182\\
2100 0.182\\
2107 0.182\\
2114 0.182\\
2121 0.183\\
2128 0.183\\
2135 0.183\\
2142 0.183\\
2149 0.183\\
2156 0.183\\
2163 0.183\\
2170 0.183\\
2177 0.183\\
2184 0.183\\
2191 0.186\\
2198 0.19\\
2205 0.194\\
2212 0.196\\
2219 0.196\\
};
\addplot [color=teal, line width=1.5pt, forget plot, line join=round, dash pattern=on 8pt off 2pt on 1.5pt off 2pt]
  table[row sep=crcr]{%
1 0.547\\
7 0.542\\
14 0.532\\
21 0.519\\
28 0.507\\
35 0.493\\
42 0.476\\
49 0.455\\
56 0.426\\
63 0.389\\
70 0.35\\
77 0.309\\
84 0.269\\
91 0.239\\
98 0.211\\
105 0.195\\
112 0.181\\
119 0.168\\
126 0.16\\
133 0.155\\
140 0.151\\
147 0.148\\
154 0.147\\
161 0.147\\
168 0.149\\
175 0.15\\
182 0.151\\
189 0.152\\
196 0.152\\
203 0.152\\
210 0.154\\
217 0.155\\
224 0.161\\
231 0.163\\
238 0.167\\
245 0.171\\
252 0.174\\
259 0.175\\
266 0.177\\
273 0.178\\
280 0.179\\
287 0.18\\
294 0.182\\
301 0.184\\
308 0.186\\
315 0.187\\
322 0.188\\
329 0.189\\
336 0.19\\
343 0.191\\
350 0.192\\
357 0.192\\
364 0.193\\
371 0.192\\
378 0.194\\
385 0.193\\
392 0.194\\
399 0.194\\
406 0.195\\
413 0.196\\
420 0.197\\
427 0.196\\
434 0.2\\
441 0.199\\
448 0.198\\
455 0.195\\
462 0.197\\
469 0.197\\
476 0.197\\
483 0.197\\
490 0.198\\
497 0.198\\
504 0.197\\
511 0.198\\
518 0.198\\
525 0.198\\
532 0.198\\
539 0.198\\
546 0.198\\
553 0.197\\
560 0.197\\
567 0.197\\
574 0.197\\
581 0.197\\
588 0.198\\
595 0.198\\
602 0.198\\
609 0.199\\
616 0.198\\
623 0.198\\
630 0.198\\
637 0.198\\
644 0.198\\
651 0.199\\
658 0.198\\
665 0.198\\
672 0.199\\
679 0.198\\
686 0.199\\
693 0.198\\
700 0.199\\
707 0.202\\
714 0.203\\
721 0.203\\
728 0.203\\
735 0.203\\
742 0.201\\
749 0.202\\
756 0.202\\
763 0.203\\
770 0.203\\
777 0.203\\
784 0.203\\
791 0.203\\
798 0.203\\
805 0.203\\
812 0.203\\
819 0.203\\
826 0.203\\
833 0.203\\
840 0.203\\
847 0.203\\
854 0.203\\
861 0.204\\
868 0.204\\
875 0.204\\
882 0.204\\
889 0.204\\
896 0.204\\
903 0.204\\
910 0.204\\
917 0.204\\
924 0.204\\
931 0.204\\
938 0.204\\
945 0.204\\
952 0.204\\
959 0.205\\
966 0.21\\
973 0.213\\
980 0.211\\
987 0.218\\
994 0.218\\
1001 0.218\\
1008 0.21\\
1015 0.205\\
1022 0.198\\
1029 0.193\\
1036 0.188\\
1043 0.184\\
1050 0.181\\
1057 0.178\\
1064 0.177\\
1071 0.176\\
1078 0.174\\
1085 0.173\\
1092 0.173\\
1099 0.171\\
1106 0.171\\
1113 0.17\\
1120 0.17\\
1127 0.169\\
1134 0.169\\
1141 0.169\\
1148 0.169\\
1155 0.168\\
1162 0.169\\
1169 0.169\\
1176 0.169\\
1183 0.17\\
1190 0.171\\
1197 0.172\\
1204 0.174\\
1211 0.177\\
1218 0.179\\
1225 0.18\\
1232 0.182\\
1239 0.183\\
1246 0.183\\
1253 0.184\\
1260 0.185\\
1267 0.185\\
1274 0.186\\
1281 0.186\\
1288 0.187\\
1295 0.187\\
1302 0.188\\
1309 0.188\\
1316 0.189\\
1323 0.189\\
1330 0.189\\
1337 0.189\\
1344 0.189\\
1351 0.189\\
1358 0.189\\
1365 0.19\\
1372 0.19\\
1379 0.19\\
1386 0.19\\
1393 0.19\\
1400 0.19\\
1407 0.191\\
1414 0.191\\
1421 0.191\\
1428 0.191\\
1435 0.191\\
1442 0.191\\
1449 0.191\\
1456 0.191\\
1463 0.191\\
1470 0.191\\
1477 0.191\\
1484 0.191\\
1491 0.191\\
1498 0.191\\
1505 0.191\\
1512 0.191\\
1519 0.191\\
1526 0.191\\
1533 0.191\\
1540 0.191\\
1547 0.191\\
1554 0.191\\
1561 0.191\\
1568 0.191\\
1575 0.192\\
1582 0.191\\
1589 0.191\\
1596 0.191\\
1603 0.191\\
1610 0.192\\
1617 0.192\\
1624 0.192\\
1631 0.192\\
1638 0.192\\
1645 0.192\\
1652 0.192\\
1659 0.192\\
1666 0.192\\
1673 0.192\\
1680 0.192\\
1687 0.192\\
1694 0.192\\
1701 0.192\\
1708 0.192\\
1715 0.192\\
1722 0.192\\
1729 0.192\\
1736 0.192\\
1743 0.192\\
1750 0.192\\
1757 0.192\\
1764 0.192\\
1771 0.192\\
1778 0.192\\
1785 0.192\\
1792 0.192\\
1799 0.192\\
1806 0.192\\
1813 0.192\\
1820 0.192\\
1827 0.192\\
1834 0.192\\
1841 0.192\\
1848 0.192\\
1855 0.192\\
1862 0.192\\
1869 0.192\\
1876 0.192\\
1883 0.192\\
1890 0.192\\
1897 0.192\\
1904 0.192\\
1911 0.192\\
1918 0.192\\
1925 0.192\\
1932 0.192\\
1939 0.192\\
1946 0.192\\
1953 0.192\\
1960 0.192\\
1967 0.192\\
1974 0.192\\
1981 0.192\\
1988 0.192\\
1995 0.192\\
2002 0.192\\
2009 0.192\\
2016 0.192\\
2023 0.192\\
2030 0.192\\
2037 0.192\\
2044 0.192\\
2051 0.192\\
2058 0.191\\
2065 0.191\\
2072 0.191\\
2079 0.191\\
2086 0.191\\
2093 0.191\\
2100 0.191\\
2107 0.191\\
2114 0.191\\
2121 0.191\\
2128 0.191\\
2135 0.191\\
2142 0.191\\
2149 0.191\\
2156 0.191\\
2163 0.191\\
2170 0.191\\
2177 0.191\\
2184 0.191\\
2191 0.191\\
2198 0.191\\
2205 0.191\\
2212 0.191\\
2219 0.191\\
};

\end{axis}

\end{tikzpicture}%
    \end{center}
    
    \textbf{Case III}\par\medskip
    \begin{center}
%
%
\definecolor{mycolor1}{rgb}{0.00000,0.44700,0.74100}%
\definecolor{mycolor2}{rgb}{0.85000,0.32500,0.09800}%
\definecolor{mycolor3}{rgb}{0.00000,0.7410000,0.32500}%
\begin{tikzpicture}

\begin{axis}[%
width=.75\linewidth,
height=.5\linewidth,
at={(0.758in,0.481in)},
scale only axis,
xmode=log,
xmin=4,
xmax=5000,
xminorticks=true,
ymode=log,
ymin=1e-3,
ymax=1,
axis background/.style={fill=white},
xmajorgrids,
xminorgrids,
ymajorgrids,
xlabel={Evaluations},
legend pos=south west
]
\addplot [color=black, line width=1.5pt, forget plot, line join=round]
  table[row sep=crcr]{%
13 1.042\\
26 0.992\\
39 0.936\\
52 0.889\\
65 0.846\\
78 0.803\\
91 0.759\\
104 0.709\\
117 0.647\\
130 0.585\\
143 0.531\\
156 0.479\\
169 0.427\\
182 0.374\\
195 0.323\\
208 0.279\\
221 0.231\\
234 0.179\\
247 0.129\\
260 0.094\\
273 0.074\\
286 0.063\\
299 0.055\\
312 0.05\\
325 0.042\\
338 0.036\\
351 0.029\\
364 0.024\\
377 0.02\\
390 0.016\\
403 0.013\\
416 0.01\\
429 0.008\\
442 0.006\\
455 0.005\\
468 0.004\\
481 0.004\\
494 0.004\\
507 0.003\\
520 0.003\\
533 0.003\\
546 0.003\\
559 0.003\\
572 0.002\\
585 0.003\\
598 0.003\\
611 0.002\\
624 0.002\\
637 0.003\\
650 0.002\\
663 0.002\\
676 0.002\\
689 0.002\\
702 0.002\\
715 0.003\\
728 0.003\\
741 0.003\\
754 0.003\\
767 0.003\\
780 0.011\\
793 0.019\\
806 0.015\\
819 0.009\\
832 0.019\\
845 0.008\\
858 0.004\\
871 0.008\\
884 0.009\\
897 0.01\\
910 0.007\\
923 0.004\\
936 0.006\\
949 0.006\\
962 0.005\\
975 0.004\\
988 0.011\\
1001 0.004\\
1014 0.006\\
1027 0.004\\
1040 0.003\\
1053 0.002\\
1066 0.002\\
1079 0.002\\
1092 0.002\\
1105 0.002\\
1118 0.002\\
1131 0.002\\
1144 0.002\\
1157 0.002\\
1170 0.002\\
1183 0.002\\
1196 0.002\\
1209 0.002\\
1222 0.003\\
1235 0.002\\
1248 0.002\\
1261 0.002\\
1274 0.002\\
1287 0.002\\
1300 0.002\\
1313 0.002\\
1326 0.006\\
1339 0.003\\
1352 0.003\\
1365 0.003\\
1378 0.003\\
1391 0.003\\
1404 0.003\\
1417 0.003\\
1430 0.003\\
1443 0.003\\
1456 0.003\\
1469 0.003\\
1482 0.003\\
1495 0.003\\
1508 0.003\\
1521 0.003\\
1534 0.003\\
1547 0.003\\
1560 0.003\\
1573 0.003\\
1586 0.003\\
1599 0.003\\
1612 0.003\\
1625 0.003\\
1638 0.003\\
1651 0.003\\
1664 0.003\\
1677 0.003\\
1690 0.003\\
1703 0.003\\
1716 0.003\\
1729 0.004\\
1742 0.004\\
1755 0.004\\
1768 0.004\\
1781 0.005\\
1794 0.006\\
1807 0.009\\
1820 0.009\\
1833 0.013\\
1846 0.036\\
1859 0.012\\
1872 0.01\\
1885 0.008\\
1898 0.007\\
1911 0.007\\
1924 0.007\\
1937 0.01\\
1950 0.007\\
1963 0.007\\
1976 0.007\\
1989 0.007\\
2002 0.007\\
2015 0.008\\
2028 0.007\\
2041 0.007\\
2054 0.007\\
2067 0.007\\
2080 0.007\\
2093 0.007\\
2106 0.007\\
2119 0.007\\
2132 0.007\\
2145 0.007\\
2158 0.007\\
2171 0.007\\
2184 0.007\\
2197 0.007\\
2210 0.007\\
2223 0.007\\
2236 0.007\\
2249 0.007\\
2262 0.007\\
2275 0.007\\
2288 0.007\\
2301 0.007\\
2314 0.007\\
2327 0.007\\
2340 0.007\\
2353 0.007\\
2366 0.007\\
2379 0.007\\
2392 0.007\\
2405 0.007\\
2418 0.007\\
2431 0.007\\
2444 0.007\\
2457 0.007\\
2470 0.007\\
2483 0.007\\
2496 0.007\\
2509 0.007\\
2522 0.007\\
2535 0.007\\
2548 0.007\\
2561 0.007\\
2574 0.007\\
2587 0.007\\
2600 0.007\\
2613 0.007\\
2626 0.009\\
2639 0.009\\
2652 0.009\\
2665 0.009\\
2678 0.009\\
2691 0.009\\
2704 0.009\\
};
\addplot [color=violet, line width=1.5pt, forget plot, line join=round, dash pattern=on 8pt off 2pt on 1.5pt off 2pt]
  table[row sep=crcr]{%
1 0.547\\
13 0.54\\
26 0.527\\
39 0.512\\
52 0.495\\
65 0.476\\
78 0.454\\
91 0.427\\
104 0.391\\
117 0.356\\
130 0.328\\
143 0.302\\
156 0.277\\
169 0.252\\
182 0.227\\
195 0.204\\
208 0.18\\
221 0.15\\
234 0.118\\
247 0.103\\
260 0.092\\
273 0.084\\
286 0.077\\
299 0.069\\
312 0.063\\
325 0.056\\
338 0.05\\
351 0.04\\
364 0.039\\
377 0.029\\
390 0.032\\
403 0.022\\
416 0.021\\
429 0.017\\
442 0.013\\
455 0.01\\
468 0.009\\
481 0.01\\
494 0.009\\
507 0.009\\
520 0.01\\
533 0.011\\
546 0.01\\
559 0.009\\
572 0.009\\
585 0.008\\
598 0.007\\
611 0.007\\
624 0.007\\
637 0.007\\
650 0.007\\
663 0.006\\
676 0.006\\
689 0.006\\
702 0.006\\
715 0.006\\
728 0.006\\
741 0.006\\
754 0.006\\
767 0.007\\
780 0.007\\
793 0.007\\
806 0.006\\
819 0.007\\
832 0.007\\
845 0.007\\
858 0.008\\
871 0.008\\
884 0.007\\
897 0.007\\
910 0.008\\
923 0.008\\
936 0.009\\
949 0.009\\
962 0.008\\
975 0.008\\
988 0.007\\
1001 0.008\\
1014 0.008\\
1027 0.007\\
1040 0.009\\
1053 0.008\\
1066 0.008\\
1079 0.008\\
1092 0.008\\
1105 0.009\\
1118 0.009\\
1131 0.009\\
1144 0.01\\
1157 0.01\\
1170 0.009\\
1183 0.009\\
1196 0.007\\
1209 0.009\\
1222 0.009\\
1235 0.009\\
1248 0.008\\
1261 0.009\\
1274 0.008\\
1287 0.008\\
1300 0.009\\
1313 0.007\\
1326 0.007\\
1339 0.007\\
1352 0.007\\
1365 0.007\\
1378 0.007\\
1391 0.008\\
1404 0.008\\
1417 0.008\\
1430 0.008\\
1443 0.008\\
1456 0.008\\
1469 0.008\\
1482 0.008\\
1495 0.008\\
1508 0.008\\
1521 0.008\\
1534 0.008\\
1547 0.008\\
1560 0.008\\
1573 0.008\\
1586 0.008\\
1599 0.008\\
1612 0.008\\
1625 0.008\\
1638 0.008\\
1651 0.008\\
1664 0.008\\
1677 0.008\\
1690 0.008\\
1703 0.007\\
1716 0.007\\
1729 0.007\\
1742 0.007\\
1755 0.007\\
1768 0.007\\
1781 0.007\\
1794 0.007\\
1807 0.007\\
1820 0.007\\
1833 0.007\\
1846 0.007\\
1859 0.007\\
1872 0.007\\
1885 0.007\\
1898 0.007\\
1911 0.007\\
1924 0.007\\
1937 0.007\\
1950 0.007\\
1963 0.007\\
1976 0.007\\
1989 0.007\\
2002 0.007\\
2015 0.007\\
2028 0.007\\
2041 0.007\\
2054 0.007\\
2067 0.007\\
2080 0.007\\
2093 0.007\\
2106 0.007\\
2119 0.007\\
2132 0.007\\
2145 0.008\\
2158 0.009\\
2171 0.007\\
2184 0.011\\
2197 0.007\\
2210 0.009\\
2223 0.007\\
2236 0.008\\
2249 0.009\\
2262 0.01\\
2275 0.007\\
2288 0.007\\
2301 0.007\\
2314 0.006\\
2327 0.006\\
2340 0.006\\
2353 0.006\\
2366 0.006\\
2379 0.006\\
2392 0.006\\
2405 0.007\\
2418 0.007\\
2431 0.007\\
2444 0.007\\
2457 0.007\\
2470 0.007\\
2483 0.007\\
2496 0.007\\
2509 0.007\\
2522 0.007\\
2535 0.007\\
2548 0.007\\
2561 0.007\\
2574 0.007\\
2587 0.007\\
2600 0.007\\
2613 0.008\\
2626 0.008\\
2639 0.008\\
2652 0.008\\
2665 0.008\\
2678 0.008\\
2691 0.008\\
2704 0.008\\
};

\end{axis}

\end{tikzpicture}%
    \end{center}
	\caption{History of convergence in all three cases. The black line indicates the approximation of the objective function~\eqref{eq:Objective} using numerical integration utilizing the evaluation of the integrand only at the targeted frequencies in each case. The dash-dotted colored lines demonstrate the true objective function approximated by numerical integration using 150 equidistant frequencies. }%
	\label{fig:MMA_History}
\end{figure}
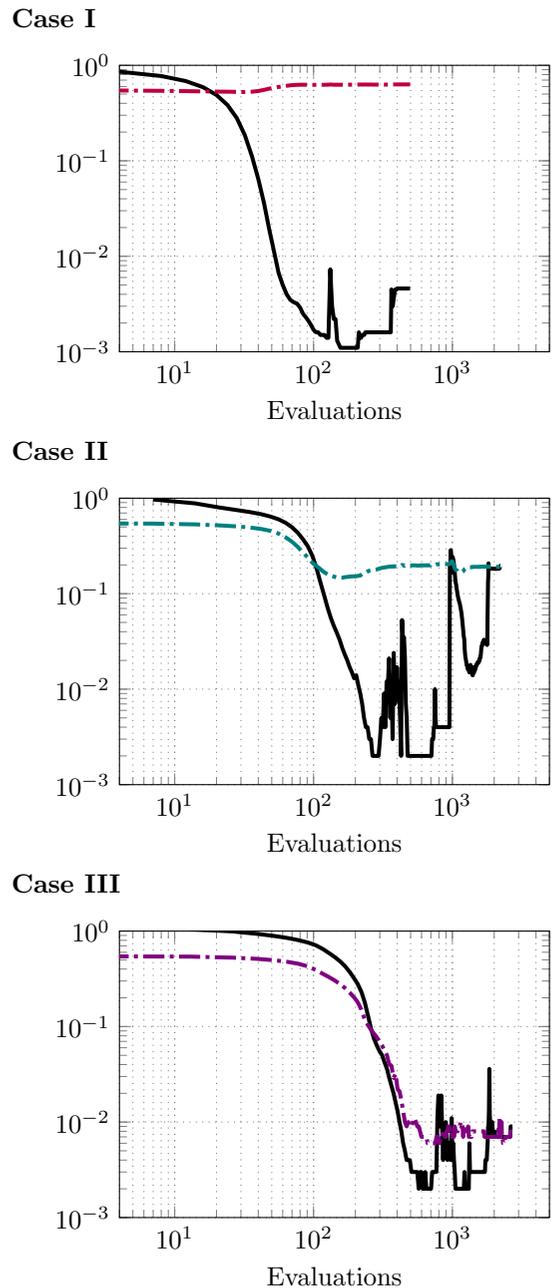

Note that here, the convergence criteria is based on the residual norm of the KKT (or the first-order optimality) condition, together with a limitation on the number of iterations in each penalty step. 
We will use the number of evaluations, defined as the number of times we need to solve the state equation~\eqref{eq:goveq}, as a metric to compare the computational cost between different cases.
Fig.~\ref{fig:MMADesigns} shows the optimized design for all three cases.
Fig.~\ref{fig:MMA} shows the performance of each of these designs, computed using expression~\eqref{eq:Performance}. 
Fig.~\ref{fig:MMA_History} shows the convergence history for all three MMA cases, where two numerical approximations of the objective function~\eqref{eq:Objective} are plotted versus the number of evaluations, one is the naive approximation using only the frequencies subject to optimization in each case and the other one is achieved using 150 equidistant frequencies in the range 4\;kHz to 16\;kHz. 

Considering the observations from Fig.~\ref{fig:MMA} and Fig.~\ref{fig:MMA_History}, the following conclusions can be drawn:
\begin{itemize}
    \item In \ref{item:Case1}, where only a few frequencies within the targeted range are considered in the optimization, the overall broadband performance of the final layout is poor, with only a few peaks observed at the considered frequencies. The total number of evaluations required for the convergence is approximately 500 evaluations.
    \item In \ref{item:Case2}, by considering more frequencies to approximate objective function~\eqref{eq:Objective}, a better overall performance is achieved for the optimized design. However, there are still deeps in the frequency response, indicating weak broadband performance. In this case, Fig.~\ref{fig:MMA_History} shows that the approximation of the objective function is still poor throughout the optimization. Moreover, approximately 2200 evaluations were required for convergence.
    \item In \ref{item:Case3}, a desirable broadband performance is achieved for the optimized layout. However, this improvement comes at the expense of an increase in computational cost, as evident from the approximately 2700 evaluations required to obtain the optimized layout. As illustrated in Fig.~\ref{fig:MMA_History}, the approximation of the objective function is significantly improved throughout the optimization compared to \ref{item:Case1} and \ref{item:Case2}.
\end{itemize}
These observations highlight two disadvantages of using a deterministic approach in this problem which can be manifested as follows:
\begin{enumerate}
    \item To ensure an acceptable broadband performance, an increased number of frequencies needs to be included in the optimization process. This results in a higher number of evaluations and computational costs.
    \item Evaluating objective function~\eqref{eq:Objective} only at specific frequencies leads to designs that are tailored for those particular frequencies. As the number of frequencies considered in the optimization increases, the resulting designs exhibit numerous free-hanging parts, as depicted in Fig.~\ref{fig:MMADesigns}. This phenomenon is characterized by an increasing number of small inclusions from \ref{item:Case1} to \ref{item:Case3} in Fig.~\ref{fig:MMADesigns}.
\end{enumerate}

In summary, while the MMA approach demonstrates advantages in achieving desirable broadband performance, it is pivotal to consider the associated computational cost and the design's specificity to the frequencies included in the optimization process.


\subsection{SG approach}
In contrast to MMA, SG does not require a computation of the integral over frequencies appearing in \eqref{eq:FinalObjective}. Instead, in each iteration, we draw a random frequency $f_n\in\mathscr{F}$ and evaluate the objective function gradient only for this specific choice. This gradient sample is then used as a search direction for the current optimization step. As a result, compared to MMA, an SG iteration consumes significantly less time, but generally provides a smaller improvement of the objective function.

While SG is typically used with a diminishing learning rate, we fix a constant learning rate and impose a shrinkage in step length directly through move limits. To be precise, in every iteration, the search direction is multiplied by the constant learning rate and afterward projected onto the set
\begin{equation*}
    \big\{ \mathbf{x}\in\mathbb{R}^{N^\text{D}}\,:\, \Vert \mathbf{x}\Vert_\infty \le C_n\big\},
\end{equation*}
where $C_n=\mathcal{O}(1/\sqrt{n})$. To end up with a black-and-white design, the final result is rounded. This setup yielded the best performance for SG in our numerical experiments.

\begin{figure}
\centering
        \input{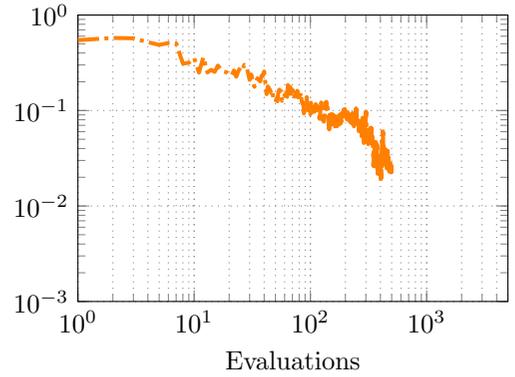}
	\caption{Value of the objective function~\eqref{eq:Objective} for all intermediate SG iterates. Each value was approximated by numerical integration using 150 equidistant frequencies.}%
	\label{fig:SG_Obj_fct}
\end{figure}

\begin{figure}
        \centering
	\includegraphics[width=\linewidth]{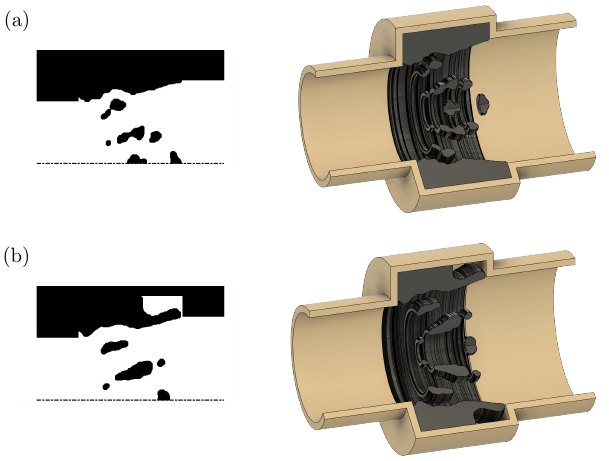}
	\caption{On the left is axisymmetric and on the right is a 3D visualization of the optimized designs in (a) SG approach and (b) CSG approach.}%
	\label{fig:SG_CSGDesigns}
\end{figure}

Due to the stochastic nature of SG, we performed 50 independent optimization runs with 500 iterations each. An overview of the observed performance range is found in Fig.~\ref{fig:500_quants}. For the best final design obtained, the objective function evolution is shown in Fig.~\ref{fig:SG_Obj_fct}, while the final design is illustrated in Fig.~\ref{fig:SG_CSGDesigns}.

\begin{figure}
\centering
        \input{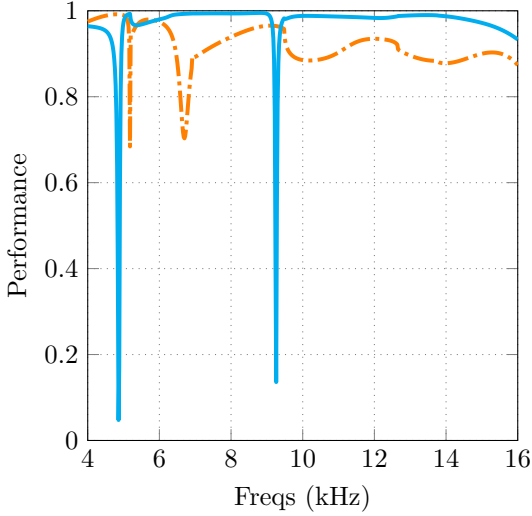}
	\caption{Performance of the optimized designs in the SG approach, the dash-dotted orange line, and in the CSG approach, the blue line. }%
	\label{fig:SG_CSG}
\end{figure}

\subsection{CSG approach}
Similar to SG, the CSG search direction is based on stochastic samples of the full objective function gradient. However, since evaluating such a gradient sample is still computationally expensive, discarding all information after each iteration is rather wasteful. Therefore, in CSG, gradient samples $(g_i)_{i=1,\ldots,n}$ from past iterations are stored. By calculating design-dependent integration weights $(\beta_i)_{i=1,\ldots,n}$, the full objective function gradient is then approximated by the continuous stochastic model
\begin{equation*}
    \nabla J(\boldsymbol{\alpha}_n) \approx \sum_{i=1}^n \beta_i g_i =: \hat{G}_n.
\end{equation*}
For our numerical experiments, we choose so-called \emph{exact hybrid weights}, since they are easily computable due to ${\operatorname{dim}(\mathscr{F})=1}$. More details on how these weights are computed in practice is given by Grieshammer et al.~\cite{CSG2}. Therein, it was also shown that the approximation error vanishes during the optimization process.
That is, 
\begin{equation*}
    \lim_{n\to\infty} \big\Vert \nabla J(\boldsymbol{\alpha}_n) - \hat{G}_n\big\Vert = 0.
\end{equation*}
As a consequence, CSG inherits strong convergence results from full gradient schemes, like convergence for constant learning rates and line search techniques, while retaining a low cost per iteration, since the integration weight computation is negligible compared to the numerical solution of the state equation.

For a better comparison, we also choose a combination of constant learning rates and move limits for CSG in our experiments. This time, however, it is not necessary to pick diminishing move limits. Instead, we can adaptively choose these limits based on the progress achieved in the internal CSG model for the objective function.

\begin{figure}
\centering
        \input{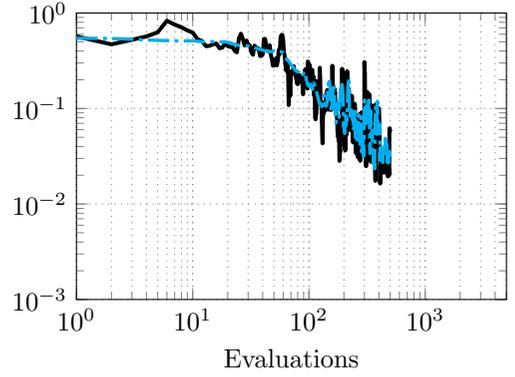}
	\caption{CSG approximation to the objective function value~\eqref{eq:Objective} in each iteration (black). The dash-dotted blue line gives a different approximation to the objective function value, obtained by using numerical integration with 150 equidistant frequencies.}%
	\label{fig:CSG_App}
\end{figure}

Again, 50 independent optimization runs with 500 evaluations each were performed. The full overview of results is shown in Fig.~\ref{fig:500_quants}. The objective function evolution of the run corresponding to the best final result obtained is given in Fig.~\ref{fig:CSG_App}. Therein, we also included the history of objective function approximations by CSG. These approximations indicate the quality of the underlying continuous--stochastic model, which is used internally in CSG. An illustration of the final design can be found in Fig.~\ref{fig:SG_CSGDesigns}.
\subsection{Impact of Stochasticity}
\begin{figure}
\centering
        \input{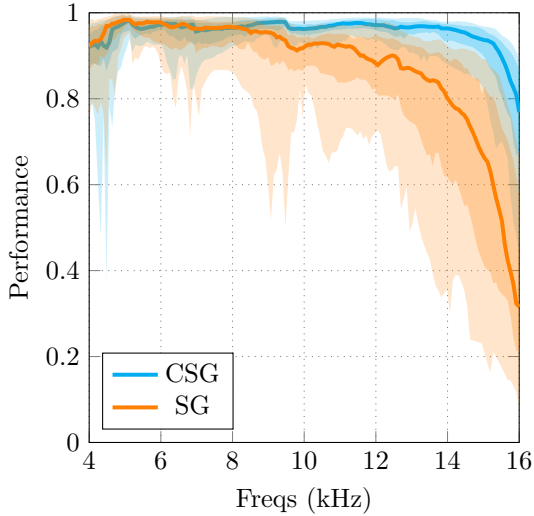}
	\caption{Median transmission spectrum (solid lines) of all final designs after 500 iterations of CSG (blue) and SG (orange). The shades areas indicate the quantiles $P_{0.1,0.9}$ (light) and $P_{0.25,0.75}$ (dark).}%
	\label{fig:500_quants}
\end{figure}

To better capture the probabilistic nature of SG and CSG, the transmission spectra of all 50 final designs are calculated. Afterwards, for each individual frequency $f\in\mathscr{F}$, we determined the respective transmission quantiles $P_{0.1,0.9}(f)$ and $P_{0.25,0.75}(f)$. Here, $P_{0.1,0.9}(f)$ denotes the range of transmission values achieved by all runs, where the highest 10\,\% and lowest 10\,\% of values are omitted. Likewise, $P_{0.25,0.75}(f)$ indicates the range of transmission values obtained by 50\,\% of all runs, where the best 25\,\% and worst 25\,\% of results are neglected. Thus, the resulting quantile plots in Fig.~\ref{fig:500_quants} give a good impression concerning both the average performance of a design obtained by SG or CSG as well as the variance in results, depending on the random sequence of sample frequencies. Note, however, that this form of representation results in smoother-looking spectra, since sharp peaks and other resonance effects are averaged out in the process.

\subsection{Discussion}
To better compare the achieved performances for all methods, we introduce the cumulative performance density $\mathrm{CPD}:[0,1]\to[0,1]$,
\begin{equation*}
    \mathrm{CPD}(p) = \frac{1}{\vert\mathscr{F}\vert}\int_{\mathscr{F}} \chi_{[0,p]}\big(P_0^\text{R}(f)\big)\,\mathrm{d}\mathscr{F}.
\end{equation*}
\begin{figure}
\centering
        \input{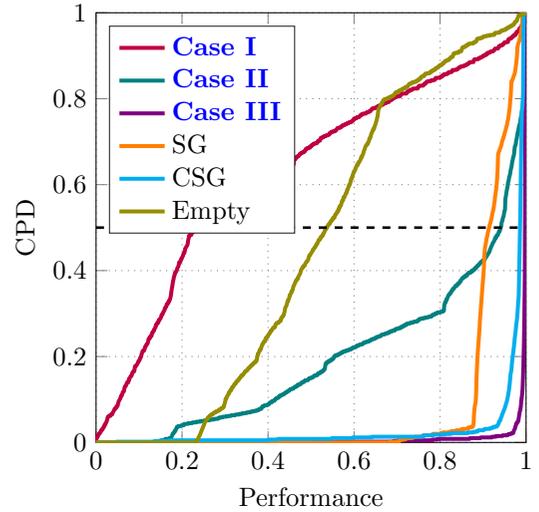}
	\caption{Cumulative performance density curves for all final designs in Fig.~\ref{fig:MMADesigns} and Fig.~\ref{fig:SG_CSGDesigns} as well as the empty design domain. For a given final design, the CPD curve indicates two measures of quality: the objective function value (area under the graph) and the median performance (intersection with horizontal line at $y=0.5$).}%
	\label{fig:CumulatedPerformances}
\end{figure}
By construction, $\mathrm{CPD}(p)$ provides the relative amount of frequencies in our frequency range $\mathscr{F}$, for which the performance is lower than the given threshold $p\in[0,1]$. For example, if a design satisfies $\mathrm{CPD}(0.5)=0.25$, its performance is less or equal to $0.5$ for $25\,\%$ of all considered frequencies. Thus, an ideal design, which has a perfect performance of $1$ for all $f\in\mathscr{F}$ satisfies
\begin{equation*}
    \mathrm{CPD}(p) = \begin{cases} 0, & p<1, \\ 1 & p=1.\end{cases}
\end{equation*}
Furthermore, the objective function value~\eqref{eq:Objective} of a design is obtained by integrating $\mathrm{CPD}$ over the interval $[0,1]$, that is,
\begin{equation*}
    J_p(\boldsymbol{\alpha}) = \int_0^1 \mathrm{CPD}(p)\dd p.
\end{equation*}

For each of the optimization approaches, the $\mathrm{CPD}$ of the corresponding final design is given in Fig.~\ref{fig:CumulatedPerformances}. Therein, we also included the $\mathrm{CPD}$ for the empty design domain ($\alpha_h \equiv 1)$. As we can see, the final design obtained with MMA in \ref{item:Case1} yields the worst performance overall, even falling behind the empty design region. The final designs of SG and MMA \ref{item:Case2} have comparable median performances. However, while the final SG design performs rather similar for all frequencies, indicated by the sharp increase of CPD at $\sim 0.9$, the MMA \ref{item:Case2} design performs poorly for a wide range of frequencies, resulting in a much better final objective function value of SG. The best overall performance is achieved by the final design of MMA \ref{item:Case3}, with the final CSG design performing slightly worse. However, recall that the associated numerical effort for MMA \ref{item:Case3} is much higher (2704 state equation solutions) when compared to CSG (500 state equation solutions).

\section{Conclusion}\label{section:CD}
In this study, we presented the results of topology optimization applied to a broadband acoustic transition section. We compared the outcomes of a deterministic approach using the method of moving asymptotes (MMA) with two stochastic approaches: stochastic gradient (SG) and continuous stochastic gradient (CSG) methods.

In the case of the MMA approach, we found that achieving  optimal broadband performance requires optimizing over an increased number of frequencies (Fig.~\ref{fig:MMA}). However, this comes at the cost of a significant increase in computational costs and results in designs with a higher number of free-hanging inclusions, which can negatively impact manufacturability (Fig.~\ref{fig:MMADesigns}).
On the other hand, the stochastic approaches, SG and CSG, offer a more computationally efficient alternative while still producing optimized designs with improved manufacturability. In particular, we observed that the CSG method outperforms the SG method, as evidenced by the median transmission spectrum of the final designs and the overall frequency response shown in the quantile plots (Fig.~\ref{fig:500_quants}).

These findings highlight the potential of stochastic approaches in acoustic applications, especially when broadband acoustic performance is desired. By reducing computational costs and improving manufacturability, stochastic methods offer a promising alternative for optimizing acoustic systems in such applications.

\section*{Declarations}
\subsection*{Funding}
A. Uihlein and L. Pflug acknowledge funding by the German Research Foundation - Project-ID 416229255 - SFB 1411.
A. Mousavi and E. Wadbro acknowledge funding by the Swedish strategic research programme eSSENCE, the Swedish Research Council (VR) (grant no.~2022-03783).

\subsection*{Conflicts of interest/Competing interests}
The authors declare that they have no known competing financial interests or personal relationships that could have appeared to influence the work reported in this paper.

\subsection*{Availability of data and material}
All the details necessary to reproduce the results have been defined in the paper.

\subsection*{Code availability}
The code is available from the corresponding author upon reasonable request.


\begin{appendices}

\section{Sensitivity analysis}\label{secA1}

Let $\boldsymbol{\alpha}$ correspond to a given design, and let $\delta \boldsymbol{\alpha}$ correspond to an arbitrary perturbation of the design.
Differentiating state problem~\eqref{eq:goveq} with respect to this perturbation, using that $\mathbf{M}$ and $\mathbf{K}$ are linear in $\boldsymbol{\alpha}$, we obtain the following expression:
\begin{multline}\label{eq:PertStat}
    \biggl(\mathbf{K}(\delta\boldsymbol{\alpha}) - k^2\mathbf{M}(\delta\boldsymbol{\alpha}) \biggr) \mathbf{p} \\
    +\biggl(\mathbf{K}(\boldsymbol{\alpha}) - k^2\mathbf{M}(\boldsymbol{\alpha}) -\mathbf{B}^\text{L} -  \mathbf{B}^\text{R}\biggr) \delta \mathbf{p} = 0.
\end{multline}
On the other hand, differentiating expression~\eqref{eq:Amplitudes_fop} yields
\begin{subequations}
\begin{align}  
    \delta B_m^\text{L} &=\lp\mathbf{v}_m^\text{L}\rp^T\mathbf{M}^\text{L} \delta \mathbf{p}, \quad m=0,\ldots,M_L^p, \label{eq:Amplitudes_fop_sen1}\\
    \delta B_n^\text{R} &= \lp\mathbf{v}_n^\text{R}\rp^T\mathbf{M}^\text{R} \delta \mathbf{p}, \quad n=0,\ldots,M_R^p. \label{eq:Amplitudes_fop_sen2}
\end{align}
\end{subequations}
Let $\mathbf{z}_m^\text{L}$ be the solution to the problem~\eqref{eq:goveq} with the right-hand side replaced by $\mathbf{M}^\text{L}\mathbf{v}_m^\text{L}$. By multiplying the linearized state problem~\eqref{eq:PertStat} with $(\mathbf{z}_m^\text{L})^T$, using the fact that the matrices are symmetric, we obtain
\begin{multline}\label{eq:Appadjoint}
    0=(\mathbf{z}_m^\text{L})^T\biggl(\mathbf{K}(\delta\boldsymbol{\alpha}) - k^2\mathbf{M}(\delta\boldsymbol{\alpha}) \biggr) \mathbf{p} \\
    +(\mathbf{z}_m^\text{L})^T \biggl(\mathbf{K}(\boldsymbol{\alpha}) - k^2\mathbf{M}(\boldsymbol{\alpha}) -\mathbf{B}^\text{L} -  \mathbf{B}^\text{R}\biggr) \delta \mathbf{p} \\
    = (\mathbf{z}_m^\text{L})^T\biggl(\mathbf{K}(\delta\boldsymbol{\alpha}) - k^2\mathbf{M}(\delta\boldsymbol{\alpha}) \biggr) \mathbf{p} \\
    + \lp\mathbf{v}_m^\text{L}\rp^T\mathbf{M}^\text{L} \delta \mathbf{p}.
\end{multline}
Substituting the last term in equation~\eqref{eq:Appadjoint} using equation~\eqref{eq:Amplitudes_fop_sen1}, we identify the partial derivatives of $B_m^\text{L}$ with respect to the element values of $\alpha_E$ in element $E$ as
\begin{equation}
    \frac{\delta B_m^\text{L}}{\delta \alpha_E} = (\mathbf{z}_m^\text{L})^T_E\biggl(k^2\mathbf{M}_E  -\mathbf{K}_E\biggr) \mathbf{p}_E,
\end{equation}
where $\mathbf{p}_E$ and $(\mathbf{z}_m^\text{L})_E$ are vectors containing the components of $\mathbf{p}$ and $(\mathbf{z}_m^\text{L})$ corresponding to nodes in element $E$ and $\mathbf{M}_E$ and $\mathbf{K}_E$ are the element mass and stiffness matrices for element $E$, respectively.

Considering $\mathbf{z}_n^\text{R}$ to be the solution to the problem~\eqref{eq:goveq} with the right-hand side replaced by $\mathbf{M}^\text{R}\mathbf{v}_n^\text{R}$, following the same argumentation as for $\mathbf{z}_m^\text{L}$, we obtain the partial derivatives of $B_n^\text{R}$ with respect to the element values of $\alpha_E$ in element $E$ as
\begin{equation}
    \frac{\delta B_n^\text{R}}{\delta \alpha_E} = (\mathbf{z}_n^\text{R})^T_E\biggl(k^2\mathbf{M}_E  -\mathbf{K}_E\biggr) \mathbf{p}_E.
\end{equation}
Note that $\mathbf{z}_m^\text{L}$ and $\mathbf{z}_n^\text{R}$ are known as adjoint variables.



\end{appendices}

\bibliography{sn-bibliography}

\end{document}